\def\a             {\alpha}
\def\Ad            {{\mathrm{Ad}}}

\def\be            {\begin{equation}}
\def\bbC           {\mathbb{C}}

\def\bbN           {\mathbb{N}}

\def\bbR           {\mathbb{R}}

\def\bbZ           {\mathbb{Z}}
\def\bfe           {{\bf1}}

\def\canr          {\theta}

\def\cA            {{\mathcal{A}}}
\def\cB            {{\mathcal{B}}}

\def\cE            {{\mathcal{E}}}

\def\cG            {{\mathcal{G}}}
\def\cH            {{\mathcal{H}}}

\def\cO            {{\mathcal{O}}}

\def\cZ            {{\mathcal{Z}}}
\newcommand\co[1]  {\overline{{#1}}}

\def\dind          {{\Delta}}

\def\E             {{\mathrm{e}}}
\def\ee            {\end{equation}}
\def\End           {{\mathrm{End}}}
\def\eps           {\varepsilon}

\newcommand\erf[1] {Eq.\ (\ref{#1})}

\def\ext           {{\mathrm{ext}}}
\def\furu          {{\mathrm{Furu}}}
\def\Gtwo          {{\mathrm{G}}_2}
\def\Hom           {{\mathrm{Hom}}}
\def\I             {{\mathrm{i}}}
\def\id            {{\mathrm{id}}}

\def\la            {\lambda}
\def\lan           {\langle}

\def\Mat           {{\mathrm{Mat}}}

\def\MXN           {{}_M {\cal X}_N}
\def\MXM           {{}_M {\cal X}_M}
\def\MXMa          {{}_M^{} {\cal X}_M^\a}
\def\MXMo          {{}_M^{} {\cal X}_M^0}

\def\MXMp          {{}_M^{} {\cal X}_M^+}
\def\MXMm          {{}_M^{} {\cal X}_M^-}
\def\MXMpm         {{}_M^{} {\cal X}_M^\pm}

\def\NXN           {{}_N {\cal X}_N}
\def\NXNd          {{}_N^{} {\cal X}_N^{\mathrm{deg}}}
\def\NXM           {{}_N {\cal X}_M}

\def\om            {\omega}
\def\Om            {\Omega}
\def\op            {{\mathrm{opp}}}

\def\ran           {\rangle}
\newcommand\res[3] {{\mathrm{res}}^{{#1}}_{{#2}}({{#3}})}
\def\rmA           {{\mathrm{A}}}
\def\rmD           {{\mathrm{D}}}
\def\rmE           {{\mathrm{E}}}

\def\RXR           {{}_R {\cal X}_R}
\def\RXRi          {{}_R^{} {\cal X}_R^\dind}
\def\RXRp          {{}_R^{} {\cal X}_R^+}
\def\RXRm          {{}_R^{} {\cal X}_R^-}
\def\RXRpm         {{}_R^{} {\cal X}_R^\pm}

\def\SLZ           {{\mathit{SL}}(2;\bbZ)}
\def\SOf           {{\mathit{SO}}(5)}

\def\SUd           {{\mathit{SU}}(3)}
\def\SUp           {{\mathit{SU}}(4)}
\def\SUf           {{\mathit{SU}}(5)}
\def\SUch          {{\mathit{SU}}(6)}

\def\SUw           {{\mathit{SU}}(8)}

\def\SUn           {{\mathit{SU}}(n)}
\def\SUN           {{\mathit{SU}}(N)}

\def\SUz           {{\mathit{SU}}(2)}

\def\SUf           {{\mathit{SU}}(4)}

\newcommand\tmat[1]{{}^{{\rm t}} {#1}}
\def\tr            {{\mathrm{tr}}}

\def\qed{{\unskip\nobreak\hfil\penalty50
\hskip2em\hbox{}\nobreak\hfil  $\Box$
\parfillskip=0pt \finalhyphendemerits=0\par}\medskip}
\def\proof{\trivlist \item[\hskip \labelsep{\it Proof.\ }]}
\def\endproof{\null\hfill\qed\endtrivlist}

\newcommand\lableq[1]{\label{#1}\end{equation}}
\newcommand\labl[1]{\label{#1}}
\def\typei         {type \nolinebreak I}
\def\typeii        {type \nolinebreak II}

\documentclass[11pt]{article}
\usepackage{amssymb,amsfonts,latexsym,epic,eepic}

\begin{document}


\newtheorem{definition}{Definition}[section]
\newtheorem{lemma}[definition]{Lemma}
\newtheorem{corollary}[definition]{Corollary}
\newtheorem{theorem}[definition]{Theorem}
\newtheorem{proposition}[definition]{Proposition}
\newtheorem{conjecture}[definition]{Conjecture}
\newtheorem{assumption}[definition]{Assumption}


\title{Critical Phenomena, Modular Invariants\\ and Operator Algebras}
\author{David E Evans\\
\\School of Mathematics\\University of Wales Cardiff\\
PO Box 926, Senghennydd Road\\Cardiff CF24 4YH, Wales, UK}
\maketitle

\begin{abstract}
We review the framework subfactors provide for understanding
modular invariants. We discuss the structure of a generalized
Longo-Rehren subfactor and the relationship between the coupling
matrices of such subfactors, modular invariance and local
extensions. We relate results of Kostant, in the context of the
McKay correspondence for finite subgroups of $\SUz$, to
subfactors.
 A direct proof of
how $\alpha$-induction produces
modular invariants is presented.
\end{abstract}

\tableofcontents

\section{Introduction}

Noncommutative operator algebras provide a framework, via the
transfer matrices, for understanding classical statistical
mechanics and their phase transitions. This has led to
connections with subfactor theory, via hyperfinite II$_1$
subfactors. Basically, algebras generated by larger and larger
partition functions with subalgebras generated by those partition
functions with certain boundary conditions. The statistical
mechanical models provide at criticality conformal field theory
models - the relevant conformal symmetry described by loop
groups, e.g those of $\SUn$ or coset models. These provide
hyperfinite III$_1$  subfactors \cite{W}. The link between the
two subfactor theories is mathematically via the classification
results of Popa \cite{Po} relating the type II$_1$ picture of
Jones-Wenzl Hecke subfactors and the type III$_1$ setting of
Jones-Wassermann loop group subfactors and physically via the
partition functions. Weyl duality, in its deformation,  relates
the representation theory of the symmetric group and the Hecke
algebras at roots of unity $q = e^{2 \pi i/(n + k)}$ to  that of
the classical group $\SUn$ quantum group $\SUn_k$ at the same
root of unity
 and the positive energy representations of
the loop group of $\SUn$ at level k. Physically, we are trying to relate
the study of the associated modular invariant
partition functions of the statistical mechanical model at criticality
and of the corresponding conformal field theories.

We start the story in more detail on the statistical mechanical
side, and for definiteness the classical Ising model on a two
dimensional square lattice $\bbZ^2$. Configuration space is
distributions of two symbols $\pm$, representing spin up/down or
occupied/unoccupied in the lattice gas picture. The Hamiltonian
$H$ is the nearest neighbour Hamiltonian. Integrability of the
model is effected mathematically by the star-triangle relation or
the Yang--Baxter equation of the Boltzmann weights
 or commutativity of the strip transfer matrices.
The transfer matrix $T$ is obtained for the partition function of
a strip of finite length $M$ and width one lattice spacing. With
boundary conditions $\xi,\eta$ along the two lengths the
corresponding partition function $T_{\xi\eta}$ defines for us the
transfer matrix $T$. The partition function $Z$ of a finite
rectangular lattice of length $M$ and width $N$ is then obtained
by multiplying the strip partition functions, namely transfer
matrix entries and summing over internal edges. For periodic
boundary conditions we obtain
\begin{equation}
\label{partition} Z=\sum_{\sigma} \exp(-\beta H
(\sigma))=\sum_{\xi_i} T{_{\xi_1\xi_2}}T{_{\xi_2\xi_3}}\ldots
T_{\xi_N\xi_1} = {\rm trace}~ T^N.
 \end{equation}

\noindent So we move from the classical model with commutative
space $\cal P$ and Ising Hamiltonian $H$ to the noncommutative
algebra, the Fermi or Pauli algebra $\bigotimes{M_2}$, and one
dimensional quantum Hamiltonian $\cH$, where $e^{- \cH} = T$ the
transfer matrix and time development $\sigma_t = Ad (T^{-it}) =
Ad (\E^{i{\cH}t)}$
 and corresponding ground states. The first step in this
movement is a noncommutative computation of the partition function
$Z = \tr T^N$ of the  $M \times N$ rectangle
(so that as $N \rightarrow \infty $, we pick up the ground states \cite{EK}).

The Ising model can be generalized to integrable models associated with
$\SUn$. The product
action of $\SUn$ on $\bigotimes{M_n}$ has a fixed point algebra generated by permutations
of the tensor factors,
or certain representation of the symmetric group. This Weyl duality and its
deformation is at the heart of the relation between the $\SUn$
loop group picture, hyperfinite III$_1$ factors, superselection
sectors and conformal
quantum field theory on the one hand and the Hecke algebra picture,
 hyperfinite
II$_1$ factors, bimodules and classicial statistical mechanics on
the other. Deforming this Weyl duality, the "fixed points" of the
quantum group action of $\SUn_q$ on $\bigotimes_{\mathbb{N}}
M_n$  is identified with a representation of the Hecke algebra
$H(q,m)$ with generators $g_i$, $i = 1, 2 \dots m$ and quadratic
relations $(q^{-1} - g_j)(q + g_j) = 0$
 and braid relations $g_i g_{i+1} g_i = g_{i+1} g_i g_{i+1}, g_ig_j =
 g_jg_i,  |i-j| > 1$.
This duality relates the labelling of the irreducible
representations of $\SUn$, (Young tableaux with fewer than $m$
rows), and the truncations to level $k$ positive energy
irreducible representations of the loop group of $\SUn$ (with the
further restriction to at most k colums) with the labelling of the
irreducible representations of the symmetric group which appear
when permuting the tensor factors, and the cut down when
$q=\E^{2\pi\I/(n+k)}$
 is a root
of unity. Let us denote either by $\cA$ =$\cA^{(n,k)}$.
The Ising Boltzman weights lie naturally in
in the Pauli algebra $\bigotimes {M_2}$
 when $q = e^{i \pi/4} $. The Yang--Baxter equation at criticality reduces to the
braid relation. The Ising model was generalized by Andrews, Baxter
and Forrester \cite{ABF}, Date, Jimbo,  Miwa, and Okado \cite{DJMO}
(to $\SUz$, $\SUn$ respectively) so that the intergrable Boltzman weights now
lie in $\bigotimes{M_n}.$ The natural labels for the states at each site of the lattice
is now $\cA^{(n,k)}$, and configuration space consists of distributions of edges
from the fusion graph $\cA^{(n,k)}$ on a square lattice $\bbZ^2$.

At criticality, the partition function $Z$ on a torus can be written as
$Z(\tau) = \tr(\E^{-\beta H}\E^{\I\eta P})$
where $2\pi\I\tau=-\beta+\I\eta$ parameterizes the metric
of the torus. The partition function of the statistical mechanical model
can be expressed as an average over $\E^{-\beta H}$,
where $\cH$  the Hamiltonian,
is now $L_0+\bar{L}_0-c/12$
(the shift by $c/24$ arising from mapping the Virasoro
algebra on the plane to a cylinder).
Here $L_0$ ($\bar{L}_0$) generate the rotation groups
which act on the underlying loop groups, and
also we have a momentum
$P$ ($=L_0-\bar{L}_0$) describing evolution along the
closed string. So taking both evolutions into account,
we compute
\be
Z(\tau) = \tr(\E^{-\beta \cH}\E^{\I\eta P}) =
\tr(\E^{2\pi \I\tau(L_0-c/24)}
\E^{-2\pi \I\bar{\tau}(\bar{L}_0-c/24)}) \,.
\labl{parf1}
\ee
Decomposing according to representation theory of the loop group symmetry,
the partition function breaks up as:
\begin{equation}
\label{Zchar}
\cZ(\tau) = \sum\nolimits_{\la,\mu} Z_{\la,\mu}
\chi_\la(\tau) \chi_\mu(\tau)^*.
\end{equation}
\noindent We want to understand the modular invariant partition function, their
realization and classification.

In the III$_1$ factor setting, we have such a factor $N$ endowed with
a braided system of endomophisms $\NXN$ - the braiding in the loop
group examples arise from rotating and possibly dilating
an interval $I$  onto its complement.
Modular invariants enter the picture when we start to look at compatable
subfactors $N \subset M$, when if $\iota$ denotes the corresponding inclusion
then the irreducible components of  the dual canonical endomorphism
 $\theta = \co \iota \iota$ lie
in $\NXN$. This setting was initated by Feng Xu \cite{X1} for
conformal embedddings where a certain chiral locality holds and taken up by B\" ockenhauer and
Evans \cite{BE1, BE2, BE3} for amongst other things simple current
or orbifold invariants and some coset models and by
 B\" ockenhauer, Evans and Kawahigashi in \cite{BE2, BE3} for so
called type II invariants when chiral locality can fail.
Xu began the ball rolling in using the method of Longo-Rehren
for inducing sectors of $N$ to those
of $M$. This is a formulation of Mackey induction, indeed the
Longo-Rehren induction reduces precisely to this in the case
of group-subgroup subfactors $R^G \subset R^H$ \cite{BE1}.
An outer finite group action on a hyperfinite
factor can be recovered from the position of the
 fixed
point algebra $R^G \subset R$. For a group-subgroup subfactor
$R^G \subset R^H$  one can not recover uniquely the subgroup $H$
\cite{KSu,I2}. Nevertheless, we think of the subfactor as
providing a mechanism for understanding the quotient $G/H$ even
when $H$ is not normal. An inclusion of a subfactor in a factor
has a rich combinatorial structure which can be described by a
quantisation of a symmetry group or a paragroup in the language
of Ocneanu \cite{O1, EK}. Jones \cite{J1} found that the minimal
algebraic structure realising this quantum symmetry was the
Temperley--Lieb algebra or quantum $\SUz$.

The induction is given by:
$\a_\la^\pm = \co\iota^{\,-1} \circ \Ad
(\eps^\pm(\lambda,\canr)) \circ \lambda \circ \co\iota$,
where $\eps$ is the braiding,
so that $\a_\la$ an endomorphism on $M$ extends $\la$ in $\NXN$
on $N$, $\a_\la\iota = \iota\la$.
The map $[\la] \rightarrow [\a_\la^\pm]$
is multiplicative, additive and preserves conjugation on sectors.
The more naive induction $\lambda \rightarrow \co \iota \lambda \iota =
{{\lambda}'}$
is not multiplicative as for example $d_{{\lambda}'}$ = $d_\theta d_\lambda$
is inconsistent with multiplicativity. We somehow need to divide
out by $\theta$ and $\alpha$-induction, in the presence of a braiding does this
so that $\a_\la \iota \co \iota = \co \iota \lambda \iota$.

We will assume initially that braiding on the system $\NXN$ is non-degenerate. In this
case there is a natural represention of the modular group
$\SLZ$ where the standard generators $S$ and $T$ matrices are basically given by the
Hopf link and twist respectively. More precisely,
recall that the statistics phase of $\om_\la$ for
$\la\in\NXN$ is given as
$ d_\la \phi_\la(\eps^+(\la,\la))=\om_\la \bfe $,
where the state $\phi_\la$ is the left inverse of $\lambda$.
We set $z =  \sum_{\la\in\NXN} d_\la^2 \omega_\la$.
If $z\neq 0$ we put $c=4\arg(z)/\pi$, which is the
central charge defined
modulo 8.
The $S$-matrix  is defined by
\[ S_{\la,\mu} =  \frac{1}{|z|} \sum_{\rho\in\NXN}
\frac{\om_\la \om_\mu}{\om_\rho} N_{\la,\mu}^\rho d_\rho \,,
\qquad \la,\mu\in\NXN \,,\]
with $N_{\la,\mu}^\rho=\lan\rho,\la\mu\ran$ denoting the
fusion coefficients  \cite{R1,FG1,FRS2}. (As usual, the
label $0$ refers to the identity morphism $\id\in\NXN$.)
Let $T$ be the diagonal matrix with entries
$T_{\la,\mu}=  \E^{-\I\pi c/12} \om_{\la} \delta_{\la\mu}$.
Then this pair of $S$ and $T$ matrices satisfy
 $TSTST=S$ and give a
unitary representation of the modular group $\SLZ$
\cite{R1, T}.
Following a study of examples in \cite{BE3}, we put
 $Z_{\la,\mu}=\lan\a^+_\la,\a^-_\mu\ran$, with the right hand
side interpreted as multiplicites of common sectors in the two inductions,
defines a matrix with positive integral entries
normalized at the vacuum, $Z_{0,0} =1$. This matrix
 commutes with $S$ and $T$ and
consequently, $Z$ gives a modular invariant
 \cite{BEK1}. (In the case of degenerate braidings, we would
define matrices $Y$ and $\Omega$, where we remove the factors,
$1/{|z|}$, and $\E^{-\I\pi c/12}$, from the above expressions for $S$ and $T$ respectively.
Then we only have the partial relations $\Omega Y\Omega Y\Omega = zY$, but
nevertheless, $\lan\a^+_\la,\a^-_\mu\ran$ still commutes with
$Y$ and $\Omega$).

A modular invariant is said to be {\it{sufferable}}
if it arises in this way from a subfactor through the process of
$\alpha$-induction. As we shall point out here again in this paper, there do exist
{\it{insufferable}} modular invariants which cannot be realised in this way.

A simple argument of Gannon \cite{G1} shows that there are
finitely many modular invariants for a given modular data. Indeed
he has the estimate $\sum_{\la,\mu}Z_{\la,\mu}\le 1/S_{0,0}^2$
which can be refined \cite{BER1} to the inequality $Z_{\la,\mu}
\le d_\la d_\mu$ for each entry of the mass matrix. There is
always the diagonal matrix solution. In some sense according to
Moore and Seiberg \cite{MS2} every modular invariant is of this
form when looked at properly. That is at least when the system
$\cal A$ is expanded to an extended system $\cal B$ and up to a
twist $\vartheta$ of the extended fusion rules of $\cal B$. So
taking the twisted diagonal modular invariant or permutation
invariant for the $\cal B$ system
\begin{equation}
Z^{\ext}=\sum\nolimits_\tau
\chi_{\tau} \co\chi_{\vartheta(\tau)},
\end{equation}
and restricting to the original system $\chi_{\tau} = \sum_{\lambda}b_{\tau,\lambda} \,\chi_\lambda $, we have:
\begin{equation}
\label{extaut}
Z_{\la,\mu}=\sum\nolimits_\tau
b_{\tau,\la} b_{\vartheta(\tau),\mu}.
\end{equation}
\noindent When no twist is necessary, so that \be \sum_{\tau\in
\mathcal{B}}|\chi_\tau|^2 = \sum_{\tau\in
\mathcal{B}}|\sum_{\lambda\in
\mathcal{A}}b_{\tau\lambda}\chi_\lambda|^2, \ee \be
Z_{\lambda\mu}~~=~~\sum_\tau b_{\tau\lambda}b_{\tau\mu}, \ee

\noindent then the modular invariant is said to be type I. In
such a case, the modular invariant is automatically symmetric
$Z_{\lambda,\mu} =
 Z_{\mu, \lambda}$. In the presence of a non trivial twist, the modular
invariant is said to be type II.
In the type II  case, the modular invariant $Z$ may not be symmetric (although
in the case of $SU(n)$ all known modular invariants are symmetric even
though not all are type I) but there is symmetric vacuum coupling:
$Z_{\lambda, 0} = Z_{0, \lambda}$. This is not enough to handle
all modular invariants -- some are not symmetric --
e.g. for the loop groups of orthogonal groups at low levels \cite{G7, BE4} or
quantum doubles of some finite groups \cite{CGR, EP}.
In this case, which one may call type III, we need different
extensions or labellings on the right and left $\cB^{\pm}$.

Such an extension theory has been developed rigorously in the
subfactor context in \cite{BE4}. At the same time a better understanding
was developed in \cite{BE4} of the fundamental formula
proposed in \cite{BE3} for $Z_{\la,\mu} = \langle \a^+_\la,\a^-_\mu \rangle$.
We write this in matrix form, $Z = b^{+*} b^-$, where $b^{\pm}$ are the
matrices of the chiral branching coefficients
\begin{equation}
b^\pm_{\tau,\la}= \lan \tau, \a^\pm_\la \ran \,,
\qquad \tau\in\cB, \quad \la\in\cA \,.
\ee
Here $\cB$ are certain induced braided neutral endomorphisms
of the factor $M$, where there are  extended $S^{ext}$ and $T^{ext}$ matrices
associated to the extended system.
Then the branching coefficient matrices $b^{\pm}$
intertwine  the two representations of $\SLZ$:
$S^\ext  b^\pm =
 b^\pm S \,, T^\ext b^\pm
=  b^\pm T$.
However, this does raise further questions, since this means
that certain modular invariants e.g. some associated
to orthogonal groups at low order are insufferable -- they cannot be realised by subfactors.
See e.g. Section 5 of \cite{BER2} for a discussion on this.
 A sufferable modular invariant with associated
subfactor produces intermediate subfactors $N\subset M_\pm\subset M$
 where neutral sectors
of endomorphims $\cB^{\pm}$ of $M_\pm$ are the natural left and right chiral
labelling producing type I parent modular invariants $Z^{\pm} =
b^{\pm*}b^{\pm}$, and the
identification of $\cB^+$ with $\cB^-$, relating $b^+$ with $b^-$,
 is the twist which yields the
original modular invariant. It may happen that $M_+\neq M_-$,
and this can even happen with symmetric modular invariants for
$\SUn$ \cite{BE4, E1}.

When we start from a braided system of endomorphisms $\NXN$ on a
III$_1$ factor $N$, and a subfactor $N \subset M$ where the dual
canonical endomorphism is generated by $\NXN$, we can induce to
systems of sectors $\MXMpm$ on $M$. We can then form the system
$\MXMa$ = $\MXMp \vee \MXMm$ and $\MXM$ as a system of
endomorphisms for the irreducible subsectors of $\{\iota \la
\co\iota : \la \in \NXN \}$ on $M$. For conformal embeddings,
where the extended sectors are given by those of the ambient loop
group, which is clearly non-degenerately braided, it was shown in
\cite{BE3}, that  the following conditions are equivalent:
\begin{itemize}
\item  the mass matrix $Z_{\la,\mu}$ can be recovered from $\alpha$-induction as
 $\langle \a^+_\la,\a^-_\mu \rangle$,

\item the neutral sectors $\MXMo$ = $\MXMp \cap \MXMm$ are identical with the
extended sectors,

\item $\MXM$ = $\MXMp \vee \MXMm$.

\end{itemize}
There were also analogous statements for simple current invariants
\cite{BE3}, and
in either case, these conditions were shown to be valid for many cases.
It was then subsequently shown in \cite{BEK1}, that if the system $\NXN$ is
non-degenerately braided, then we have the generating property
$\MXMa$ = $\MXM$ (and hence the equivalent conditions of
\cite{BE3} hold). A more algebraic proof of the generating property
which does not rely on graphical
arguments and the double triangle algebra was established in \cite{BE4}. In the same spirit,
we give an algebraic proof in \ Section \ref{algebraic}, that the mass matrix $Z$
commutes with the $S$ matrix, and is consequently a modular invariant
as a simple argument of \cite[Lemma 6.1]{BE4}) already gives commutativity with $T$.
Now in general, the induced systems, $\MXMpm$ are not braided or even
commutative. Indeed the complexifications of the finite dimensional
C*-algebras they generate are given as \cite{BEK2}
using the branching coefficients:

\begin{equation}
\label{chirfudec} \furu(\MXMpm) \simeq \bigoplus_{\la\in\NXN}
\bigoplus_{\tau\in\MXMo} \Mat(b^\pm_{\tau,\la}).
\end{equation}

In particular, the induced systems are noncommutative
only when the branching coefficients are all at most $1$.
The first noncommutative example was computed by hand
by Feng Xu \cite{X1}, for the conformal embedding
$SU(4)_4 \subset SU(15)_1$, but we see that this is
just the first of a series of noncommutative examples
from $SU(n)_n \subset SU(n^2 -1)_1$. So naturally,
using $\alpha$-induction, we have an action or nimrep
(i.e.\ a matrix representation where all the matrix entries
are non-negative integers) of the
original fusion rules of $\NXN$ on $\MXMpm$. However, this does not
yield the correct graph in general. Computing dimensions in \erf{chirfudec}
we see that the number of irreducible components of
$\MXMpm$ is given by $\tr Z^{\pm}$ - which only sees the type I parents of
the original invariant. Consider the three modular invariants
associated with the Dynkin diagrams $\mathrm{A}_{17}$, $\mathrm{D}_{10}$ and $\mathrm{E}_7$ for
$\SUz$ at level 16:
\[
\begin{array}{lll}
\cZ_{\mathrm{A}_{17}} &= \sum\nolimits_{\la=0}^{16} |\chi_\la|^2\\
\cZ_{\mathrm{D}_{10}} &=  |\chi_0 + \chi_{16}|^2
+ |\chi_2 + \chi_{14}|^2 + |\chi_4 + \chi_{12}|^2
+ |\chi_6 + \chi_{10}|^2 + 2|\chi_8|^2,\\
\cZ_{\mathrm{E}_7} &= |\chi_0 + \chi_{16}|^2
+ |\chi_4 + \chi_{12}|^2 + |\chi_6 + \chi_{10}|^2
+ |\chi_8|^2 \\
&\qquad\qquad + (\chi_2 + \chi_{14})\chi_8^*
+ \chi_8 (\chi_2 + \chi_{14})^*.
\end{array}
\]
For either $\cZ_{\mathrm{D}_{10}}$ or  $\cZ_{\mathrm{E}_7}$
 we would only get 10 sectors which would
force $\rmD_{10}$ to be the fusion graphs of $\MXMpm$.
To distinguish the two invariants, we need to
consider not only $N$-$N$ sectors $\NXN$ and $M$-$M$ sectors $\MXM$ but
also $M$-$N$ sectors $\MXN$ given by irreducible components of
$\{\iota\la: \la \in \NXN \}$.
The left action of $\MXM$ on $\MXN$ defines a
representation $\varrho$ of the $M$-$M$ fusion rule algebra,
and decomposes as
\cite{BEK1,BEK2}
\[
\varrho \simeq \bigoplus_{\la\in\NXN} \pi_{\la,\la} \,.
\]
Here $\pi_{\la,\la}$ is the irreducible representation
corresponding to the matrix block $\Mat(Z_{\la,\la})$,
$\pi_{\la,\la}([\a^\pm_\nu])=S_{\la,\nu}/S_{\la,0}\bfe_{Z_{\la,\la}}$
in the diagonalisation of the complexification of the fusion rule
algebra of $\MXMpm$
(also provided that the braiding on $\NXN$ is non-degenerate):
\begin{equation}
\label{fullfudec}
\furu(\MXM) \simeq \bigoplus_{\la,\mu\in\NXN}
\Mat(Z_{\la,\mu}).
\end{equation}
In particular the spectrum is determined by the diagonal part
of the modular invariant. Thus it is precisely this representation
$\varrho$ which provides an automatic connection between
the modular invariant and fusion graphs (e.g.\ the
representation matrix of some fundamental generator
$\Box$ corresponding to the left multiplication
of $[\a^\pm_{\Box}]$ on the $M$-$N$ sectors)
in such a way that (the multiplicities in) their spectra are
canonically given by the diagonal entries of the coupling matrix.
In fact, evaluation of $\varrho$ on the $[\a^\pm_\la]$'s
yields a nimrep of the original $N$-$N$ fusion rules.

We say a modular invariant is {\it{nimble}} if there exists
a nimrep $\{G_\lambda : \lambda \in \NXN \}$ whose spectrum
reproduces the diagonal part of the modular invariant
$\sigma(G_\lambda) = \{S_{\lambda\mu}/S_{\mu 0}$ with multiplicity
$Z_{\mu\mu} \}$. A sufferable modular invariant is clearly
nimble by \cite{BEK1}, but not all modular invariants are nimble
and there do exist nimble insufferable modular invariants \cite{EP}.

Finally, by counting dimensions we see that $\#\MXN=\tr(Z)$.
Thus for $\SUz$ at level 16 we have
 $\tr \cZ_{\rmD_{10}} = 10$, $\tr \cZ_{\rmE_7}  = 7$,
 so that we would get the
correct $\rmE_7$ graph for the $\rmE_7$ modular invariant here. For
the situation where chiral locality holds as for $\rmD_{10}$,
we can identify both $\MXMpm$ with $\NXM$ (by
$\beta\mapsto\beta\circ\iota$, $\beta\in\MXMpm$)
so that there is no difference
here with the first cruder computation.

It may be useful here to list the exceptional invariants
for $\SUz$, $\SUd$, $\SUp$ - i.e. the ones which are not diagonal,
orbifold or simple current invariants or their conjugates.
If $Z$ is a modular invariant, then $Z^c$ = $ZC$ = $CZ$
is also a modular invariant as $C = S^2$ ,
called the conjugate of $Z$ which may different from $Z$.
Also, first for the record all $\SUn_k$ conformal embeddings
are:
\begin{eqnarray*}
&\SUn_{n-2}\subset {{\mathit{SU}}(n(n-1)/2)}_1\ ({\rm only\
gives\ exceptionals\ for\ }n\ge 5)\,;\\
&\SUn_n\subset {{\mathit{SO}}(n^2-1)}_1\ ({\rm only\ gives\ exceptionals\
for\ }n\ge 4)\,;\\
&\SUn_{n+2}\subset {{\mathit{SU}}(n(n+1)/2)}_1\ ({\rm only\ gives\
exceptionals\ for\ }n\ge 3)\,;\\
&\SUz_{10}\subset \SOf_1\,;\quad
\SUz_{28}\subset (\Gtwo)_1\,;\\
&\SUd_9\subset (\rmE_6)_1\,;\quad
\SUd_{21}\subset (\rmE_7)_1\,;\\
&\SUp_8\subset {{\mathit{SO}}(20)}_1\,;\quad
\SUch_6\subset {{\mathit{Sp}}(20)}_1\,;\quad
\SUw_{10}\subset {{\mathit{SO}}(70)}_1\,.
\end{eqnarray*}

\medskip

\noindent $\SUz$

\begin{itemize}

\item $\rmE_6$,   $\SUz_{10}\subset\SOf_1$

\item $\rmE_8$,   $\SUz_{28}\subset(\Gtwo)_1$

\item $\rmE_7$, automorphism or twist of the orbifold invariant $\rmD_{10}$ =
$\SUz_{16}/\bbZ_2$

\end{itemize}

\noindent $\SUd$

\begin{itemize}

\item  $\SUd_3\subset\mathit{SO}(8)_1$, also realised as an
 orbifold $\SUd_{3}/\bbZ_3$

\item ${\cE^{(8)}}$, $\SUd_5\subset\mathit{SU}(6)_1$, plus conjugate
${\cE^{(8)c}}$ =${\cE^{(8)}}/\bbZ_5$

\item ${\cE^{(12)}}$, $\SUd_9\subset (\mathit{E}_6)_1$, with two nimreps,
${\cE^{(12)}}$, and ${\cE^{(12)}}/\bbZ_3$,

\item ${\cE^{(24)}}$, $\SUd_{21}\subset(\mathit{E}_7)_1$

\item  ${\cE^{(12)}_{MS}}$ Moore-Seiberg invariant, automorphisim or twist of the  orbifold
invariant $\SUd_{9}/\bbZ_3$, plus conjugate ${\cE^{(12)c}_{MS}}$

\end{itemize}

\medskip

\noindent $\SUp$

\begin{itemize}

\item ${\cE^{(8)}}$ $\SUp_4\subset\mathit{SO}(15)_1$

\item ${\cE^{(10)}}$ $\SUp_6\subset\mathit{SU}(10)_1$, plus conjugate
${\cE^{(10)c}}$ = ${\cE^{(10)}}/\bbZ_5$

\item ${\cE^{(12)}}$ $\SUp_8\subset\mathit{SO}(20)_1$

\item ${\cE^{(12)}_{aut}}$ $\SUp_8$: automorphism or twist of the  orbifold invariant
$\SUp_{8}/\bbZ_4$

\end{itemize}

Of course apart from trying to realise modular invariants by subfactors
or vice versa one would like to classify all modular invariants for a
given modular data, and understand the relevant subfactors  which
can be used to realise them through the process of $\alpha$-induction.

Cappelli, Itzykson and Zuber \cite{CIZ2} classified $\SUz$ modular
invariants by ADE Dynkin diagrams whose exponents or eigenvalues
describe the diagonal part of the invariant. Such graphs are
obtained from the McKay graphs for finite subgroups of $\SUz$,
the affine Dynkin diagrams. Di Francesco and Zuber \cite{DZ1}
understood the ADE classification of $\SUz$ modular invariants
\cite{CIZ2} in terms of nimreps. They were very successful in
pursuing this line in trying  to label the known $\SUd$ modular
invariants (subsequently classified by Gannon \cite{G2}) with
graphs whose eigenvalues described the diagonal part of the
modular invariant, partly guided by the list of McKay fusion
graphs for finite subgroups of $\SUd$. There is a mismatch
between the lists of $\SUd$ modular invariants and finite
subgroups of $\SUd$ but they are very close and this
"correspondence" was a significant tool.
 As $n$ increases, then the number of finite subgroups of
$\SUn$ increases,
but the corresponding number of exceptional modular invariants appears to
decrease.
They succeeded, for $\SUd$, in finding graphs and nimreps for the orbifold invariants,
and the exceptional invariants, (with three candiates for the
conformal embedding  $\SUd_9 \subset (\mathit{E}_6)_1$
invariant). All these graphs were three-colourable, and they conjectured
this to be the case for all $\SUd$ modular invariants.
Meanwhile, an ADE classification of subfactors of index less than
four was obtained by Ocneanu \cite{O1} with published proofs by
others (see \cite{EK} for detailed citations).
Here $\mathit{D}_{odd}$ and $\mathit{E}_7$ are missing from this
classification as they do not appear as the principal
graphs of any subfactor. These correspond to the {\typeii}  modular invariants in the
corresponding ADE classification. The orbifold construction was brought
into the subfactor picture in \cite{EK1}, and followed up in \cite{Kaw, X3}.
It should be noticed that the
flat part of the  $\mathit{E}_7$ connection was computed to be
 $\mathit{D}_{10}$ in
\cite{EK2}, motivated by the fact that
the $\mathit{E}_7$  modular invariant was a twist on the orbifold
$\mathit{D}_{10}$  invariant. These were pieces of a jigsaw
relating the classification of $\SUn$ modular invariants with
$\SUn$ subfactors, but the overall picture was fragmented.

Looking further into the subfactor connections,
Ocneanu \cite{O} produced the nimreps, chiral and full graphs for the $\SUz$ invariants
using the bimodule and his double triangle approach.
Feng Xu \cite{X1} looked at the conformal embedding invariants in  the loop group setting
of \cite{W}, taking as principal tool the
$\a$-induction of Longo-Rehren \cite{LR}.
In particular he computed the chiral graphs for the $\rmE_6$,
$\rmE_8$ for $\SUz$ and  $\SUd_3\subset\mathit{SO}(8)_1$,
${\cE^{(8)}}$,  ${\cE^{(12)}}$,
 ${\cE^{(24)}}$ for $\SUd$, and ${\cE^{(8)}}$,
 $\SUp_4\subset\mathit{SO}(15)_1$ for $\SUp$, and $\SUf_3 \subset
\mathit{SU}(10)_1$.
B\"ockenhauer and Evans \cite{BE1, BE2, BE3}
pursuing the $\alpha$-
 induction approach also brought the simple current invariants into the
 game and so computed the chiral graphs and full graphs in the orbifold cases
 $\mathit{D}_{even}$, as well as the
 full graphs in the case of  $\rmE_6$,
$\rmE_8$. It only then became clear that the classification
obtained by Ocneanu of irreducible connections on the ADE Dynkin
diagrams which form a fusion ring with connection generators $W$ and $\co W$
was related to the induced full $M$-$M$ system with generators
$\alpha^+_1$ and $\alpha^-_1$. The precise relation between the
bimodule approach of Ocneanu and that of $\alpha$-induction is
explained mathematically in \cite{BEK1}. The principal and dual
principal graphs could be computed from the induced systems, e.g.
 the dual principal graph for the conformal
embedding  ${\cE^{(8)}}$: $\SUd_5\subset\mathit{SU}(6)_1$
was computed for the first time  \cite{BE3}.
 The centre $\bbZ_n$
of $\SUn$ acts on the
algebra $N = \pi_0({\mathit{L}}_I{\SUn})''$, for say the vacuum
level $k$ representation, for the $\SUn$ valued
loops concentrated on an interval $I$ on a circle.
 We can form
the crossed product subfactor $N(I) \subset N(I) \rtimes \bbZ_n$,
which will recover the orbifold modular
invariants, but this extended system is only local
 if and only if $k\in 2n\bbN$
if $n$ is even and $k\in n\bbN$ if $n$ is odd \cite{BE2},
precisely when the corresponding modular invariants are \typei.
They also  started to consider the
modular invariants of the extended $U(1)$ current algebras
as treated in \cite{bumt} and the minimal model
modular invariants which arise from coset theories
$(\SUz_{m-2}\otimes\SUz_1)/\SUz_{m-1}$ and are labelled by
pairs $(\cG_1,\cG_2)$ of ADE graphs, associated to
levels $(m-2,m-1)$, which arise from fusion graphs from
$\a$-induction. See \cite{KL} for the latest developments.

The situation when chiral locality \cite{BE3} fails was taken
by B\"ockenhauer, Evans and Kawahigashi \cite{BEK2}. In particular they
realised the $\mathit{E}_7$ invariant for $\SUz$ through a subfactor,
computed the full system, as well as the full system for
$\mathit{D}_{odd}$. Previously, at a purely algebraic level,
B\"ockenhauer and  Evans \cite{E2} had constructed the nimrep for
the $\mathit{E}_7$ type Moore Seiberg invariant, a twist
of an orbifold invariant,  for $\SUd$
without realising the subfactor - but it is the sufferablity
of such invariants which is of main interest not simply constructing
a nimrep.
B\"ockenhauer and  Evans \cite{B3}
understood that nimrep graphs for the conjugate
$\SUz$ modular invariants were not three colourable.
This was also realised simultaneously by Behrend, Pearce, Petkova and
Zuber \cite{BPPZ} and Ocneanu \cite{O2}.
Indeed Ocneanu announced in Bariloche \cite{O2}
that all $\SUd$ modular invariants are sufferable, and
the classification of their associated nimreps.
B\"ockenhauer and  Evans \cite{BER1}
realised all modular invariants for cyclic
$\bbZ_n$ theories, in particular charge conjugation.
This could be used to understand and realise
 $\SUd_9\subset (\mathit{E}_6)_1$, with two nimreps.
One was ${\cE^{(12)}}$  through of course
the  $\SUd_9\subset (\mathit{E}_6)_1$
standard conformal embedding,
and another the orbifold ${\cE^{(12)}}/\bbZ_3$ obtained from the
subfactor $\SUd_9\subset (\mathit{E}_6)_1 \rtimes \bbZ_3$.
The extension $(\mathit{E}_6)_1 \subset (\mathit{E}_6)_1 \rtimes \bbZ_3$
describing charge conjugation on the cyclic $\bbZ_6$ system for
$(\mathit{E}_6)_1$. The analogous
construction for $\SUp$ does not degenerate in this way as
the modular invariant ${\cE^{(10)}}$:  $\SUp_6\subset\mathit{SU}(10)_1$,
is not self conjugate, and so
the fusion graph for the conjugate modular invariant
is the $\bbZ_5$ orbifold \cite{BER1}. A nimrep for the
conformal embedding
 invariant ${\cE^{(10)}}$: $\SUp_6\subset\mathit{SU}(10)_1$ had been
computed  in \cite{PZ2}.
For other reviews
and further  results  see \cite{BER1, BER2}.

On the physical side the role of nimreps became clearer through
the work first of Behrend, Pearce, Petkova and Zuber \cite{BPPZ}
on boundary conformal field theory bCFT and later Petkova and
Zuber \cite{PZ3, PZ4},
 particularly in the
use of twisted partition functions and its relation to the
Ocneanu approach. On the other hand the $\alpha$-induction approach
or formulation was taken up on the physical side by F\"uchs
and Schweigert \cite{FS2, FS3} and adapted to their approach to bCFT.

Returning to our decomposition of modular invariants through
intermediate subfactors, the
situation is summarised \cite{BE3} using recent work of Rehren
\cite{R7} on canonical tensor product subfactors as
determining a generalized Longo Rehren LR subfactor
with a maximal intermediate extension,
\be
Q=N\otimes N^\op \subset M_+\otimes M_-^\op \subset R,
\label{inter}
\ee
where the irreducible sector decomposition of the
dual canonical endomorphism $\Theta$ of $Q\subset R$
is described by $Z$ \erf{redce}, and $\Theta^\ext$ of
$M_+\otimes M_-^\op\subset R$ by $Z^\ext$.

 There is a  connection between the two chiral
inductions and the picture of left- and right-chiral algebras in
conformal field theory. Suppose that our factor $N$ is obtained as
a local factor $N=N(I_\circ)$ of a quantum field theoretical net
of factors $\{N(I)\}$ indexed by proper intervals $I\subset \bbR$
on the real line, and that the system $\NXN$ is obtained as
restrictions of DHR-morphisms (cf.\ \cite{H}) to $N$. This is in
fact the case in our examples arising from conformal field theory
where the net is defined in terms of local loop groups in the
vacuum representation. Taking two copies of such a net and placing
the real axes on the light cone, then this defines a local net
$\{A(\cO)\}$, indexed by double cones $\cO$ on two-dimensional
Minkowski space (cf.\ \cite{R5} for such constructions). Given a
subfactor $N\subset M$, determining in turn two subfactors
$N\subset M_\pm$ obeying chiral locality, will provide two local
nets of subfactors $\{N(I)\subset M_\pm(I) \}$ as a local
subfactor basically encodes the entire information about the net
of subfactors \cite{LR}. Arranging $M_+(I)$ and $M_-(J)$ on the
two light cone axes defines a local net of subfactors
$\{A(\cO)\subset A_\ext(\cO)\}$ in Minkowski space. The embedding
$M_+\otimes M_-^\op \subset B$ gives rise to another net of
subfactors $\{A_\ext(\cO) \subset B(\cO)\}$, where the net
$\{B(\cO)\}$ obeys local commutation relations. The existence of
the local net was already proven in \cite{R7}, and now the
decomposition of $[\Theta_\ext]$ tells us that the chiral
extensions $N(I)\subset M_+(I)$ and $N(I)\subset M_-(I)$ for left
and right chiral nets are indeed maximal (in the sense of
\cite{R5}), following from the fact that the coupling matrix for
$\{A_\ext(\cO) \subset B(\cO)\}$ is a bijection. This shows that
the inclusions $N\subset M_\pm$ should in fact be regarded as the
subfactor version of left- and right maximal extensions of the
chiral algebra.

The extension \erf{inter}
 basically ensures the existence of a corresponding
two dimensional rational conformal field theory 2D-RCFT \cite{R7}.
The coupling matrix $Z$ is automatically $\Om$- and
$Y$-invariant \cite{BEK1}.
Crucial for the understanding of classifications of
2D-RCFT's and braided subfactors is now the converse
direction: If we have generalized LR-subfactor $Q\subset R$
(relative to a braided system), i.e.\ a 2D-RCFT,
does this also imply the existence of some braided
subfactor $N\subset M$ such that its $\a$-induction
coupling matrix reproduces the $\Theta$ of $Q\subset R$?
A proper answer to this question would clarify
why for RCFT's the diagonals are exponents of graphs etc.
This question is not easy. The first problem is that
we cannot expect a unique subfactor $N\subset M$ because
different subfactors can produce the same $Z$, see
\cite{BEK2,BER1,BER2} and the changing the $\iota$-vertex
argument \cite{E1}. The next problem is that even
a generalized LR subfactor which satisfies the 2D locality
condition produces a coupling matrix which is $\Om$-invariant
but which does not necessarily commute with $Y$.
The most trivial example is the trivial (and hence local)
inclusion $Q=R$ which gives
$Z_{\la,\mu}=\delta_{\la,0}\delta_{\mu,0}$ which is
obviously $\Om$- but not $Y$-invariant.
(Rehren has  a less trivial example for $\SUz_k$:
There is a local $Q\subset R$ such that the only
non-zero entries of $Z$ are $=1$ on the diagonal
for even spins. The locality condition is possibly
equivalent to but at least should imply $\Om$-invariance of the coupling matrix.)
However, it follows from Proposition \ref{locint} that the
inclusions $N\subset M_\pm$ determined by $Q\subset R$ are
braided subfactors and hence produce $\Om$- \textit{and}
$Y$-invariant (parent) coupling matrices but which can
(and often are) different from the coupling matrix
describing $\Theta$. (These are however trivial as
$N=M_\pm$ for the sketched examples.)
So the question is: When does a generalized LR subfactor
$Q\subset R$ with coupling matrix $Z$ determine some
braided subfactor $N\subset M$ giving $Z$ by $\a$-induction?
Is it enough to require $Y$-invariance? (Besides locality
ensuring $\Om$-invariance.)

Clearly Lie groups, and their loop groups, provide
clear examples for studing modular data and modular invariants.
However, finite discrete groups as we have already mentioned,
have a role to play in this theory - from finite subgroups at
least of $\SUz$ and $\SUd$. This aspect will be further discussed
in \ Section \ref{McKay}
when we relate the approach of Kostant \cite{Kos}
to our subfactor approach. Apart from loop group or quantum group
examples, modular data can also be obtained in the subfactor context
through the use of the quantum double subfactor, just as one
can produce $R$-matrices from any Hopf algebra through the
double construction of Drinfeld. Thus the quantum doubles of
finite groups provide a fascinating testing ground \cite{EP} for testing
whether structure or coincidences which we have discovered
in the quantum group subfactors, carry over to more general
braided subfactors and their corresponding modular invariants.

\section{On the structure of a generalized Longo-Rehren subfactor}
\labl{GLR}

Here we collect some considerations about the relation
between braided subfactors and generalized Longo-Rehren
(LR) subfactors. Roughly speaking, the trivial subfactor
$N\subset N$, as producing the diagonal coupling matrix,
corresponds to the usual LR subfactor. More interesting
are non-trivial inclusions $N\subset M$ producing
non-diagonal coupling matrices and corresponding to
\textit{generalized} LR subfactors.

\subsection{Maximal local extensions}

The following  will demonstrate that a
generalized LR subfactor gives naturally
the local extensions factor $M_\pm$.

Let $N$ be a type III factor with a system $\NXN$ of
endomorphisms, and let $Q=N\otimes N^\op$.
We consider a subfactor $Q\subset R$ such that the
corresponding dual canonical endomorphism is
decomposed as
\be
[\Theta] = \bigoplus_{\la,\mu\in\NXN} Z_{\la,\mu}
[\la\otimes\mu^\op],
\label{redce}
\ee
i.e.\ $Q\subset R$ is a canonical tensor product subfactor
\cite{R5,R7}. Let $\iota_\mathrm{LR}:Q\hookrightarrow R$ be the
injection homomorphism, so that
$\Theta=\co\iota_\mathrm{LR}\iota_\mathrm{LR}$, and choose
orthonormal bases of isometries
\[
\Psi_{\la,\mu;i} \in \Hom(\iota,\iota\circ(\la\otimes\mu^\op)),
\qquad i=1,2,\ldots,Z_{\la,\mu} \,.
\]
Note that then any $r\in R$ has a unique expansion in the
$\Psi_{\la,\mu;i}$'s with coefficients in $Q$ \cite{LR,ILP}.
In the sequel we will freely identify $N$ canonically with
$N\otimes\bfe$ and similarly $N^\op$ with $\bfe\otimes N^\op$.

\begin{proposition}
\label{locint}
There are irreducible subfactors $N\subset M_\pm$ such that
their dual canonical endomorphisms $\canr_\pm$ are decomposed
as $[\canr_+]=\bigoplus_{\la\in\NXN} Z_{\la,0}[\la]$.
and $[\canr_-]=\bigoplus_{\la\in\NXN} Z_{0,\la}[\la]$,
respectively.
\end{proposition}

\proof
Put $\hat N = N \otimes \bfe \subset R$ and define
\[
\hat{M}_+ = \sum_{\la\in\NXN} \sum_{i=1}^{Z_{\la,0}}
\hat N \Psi_{\la,\mu;i}.
\]
We first claim that $\hat N \subset \hat{M}_+$ is
an irreducible subfactor, i.e.\
$\hat N ' \cap \hat M _+= \bbC\bfe$
(in particular $\hat M _+$ is a factor).
Let $X\in\hat N ' \cap \hat M _+$ and expand
$X=\sum_{\la,i} (n_{\la,i}\otimes\bfe)\Psi_{\la,0;i}$,
with $n_{\la,i}\in N$.
Since the $\Psi_{\la,\mu;i}$'s are linearly independent
over $Q$ we obtain from $X\in \hat N '$ for all $n\in N$
\[
(n\otimes\bfe) (n_{\la,i}\otimes\bfe) \Psi_{\la,0;i}
= (n_{\la,i}\otimes\bfe) \Psi_{\la,0;i} (n\otimes\bfe)
= (n_{\la,i}\otimes\bfe) (\la(n)\otimes\bfe) \Psi_{\la,0;i}
\]
for all $\la\in\NXN$, $i=1,2,...,Z_{\la,0}$ individually,
where we also used the intertwining property of
$\Psi_{\la,0;i}$. We conclude
$n_{\la,i}\in\Hom(\la,\id)=\delta_{\la,0}\bbC\bfe$,
implying $\hat N ' \cap \hat M _+= \bbC\bfe$.
Letting $\iota_+:\hat N \hookrightarrow \hat M _+$ denote
the injection homomorphism we find with a similar argument
for $\la'\in\End(N)$ that
$\Hom(\iota_+,\iota_+(\la'\otimes\id))$
is zero unless $[\la']=[\la]$ and consequently we obtain
by \cite{ILP} that $\canr_+= \iota_+\iota_+$ is
decomposed as
$[\canr]=\bigoplus_\la Z_{\la,0} [\la\otimes\id]$.
Playing a similar game with
$\hat M _-^\op = \sum_{\mu,i} (\bfe\otimes N^\op) \Psi_{0,\mu;i}$
completes the proof.
\endproof

The following identifies $\hat M _+$ and $\hat M _-^\op$
as relative commutants of $N^\op$ and $N$, respectively, in $R$,
and therefore can be considered as some kinds of
``maximally extended chiral algebras'' (cf.\ \cite{R5}).

\begin{proposition}
We have
\be
\hat M _+ = (\bfe\otimes N^\op)' \cap R , \qquad
\hat M _-^\op = (N\otimes \bfe)' \cap R .
\ee
\end{proposition}

\proof
We only show the first equality, the second one is proven similarly.
Any $X\in R$ can be expanded as
\[
X = \sum_{\la,\mu,i\in\NXN} \sum_{i=1}^{Z_{\la,\mu}}
(n^1_{\la,\mu;i} \otimes n^2_{\la,\mu;i}) \Psi_{\la,\mu;i}
\]
with $n^1_{\la,\mu;i}\in N$, $n^2_{\la,\mu;i}\in N^\op$.
Then $X\in (\bfe\otimes N^\op)'$ implies that
$n^2_{\la,\mu;i}\in\Hom(\mu^\op,\id^\op)=\delta_{\mu,0}\bbC\bfe$,
so that $X\in\hat M _+$. Consequently,
$(\bfe\otimes N^\op)' \cap R \subset \hat M _+$.
The opposite inclusion is obvious.
\endproof

>From now on we assume that the system $\NXN$ is braided
and denote
\[
\hat{\eps}^\pm(\la\otimes\mu^\op,\rho\otimes\nu^\op)^* =
\eps^\pm(\la,\rho)\otimes j(\eps^\pm(\mu,\nu)).
\]
(Here $j:N\rightarrow N^\op$ denotes the anti-linear
isomorphism.)
This gives a braiding on $\NXN\otimes\NXN^\op$ and
can be extended to sums and products in the usual way
(see e.g.\ \cite[Subsect.\ 2.2]{BEK1}).

\begin{proposition}
If the locality condition holds,
$\hat{\eps}^\pm(\Theta,\Theta)V^2=V^2$ for
Longo's isometry $V\in\Hom(\id,\Gamma)$,
$\Gamma=\iota\co\iota$, then $\hat M _+$
and $\hat M _-^\op$ commute.
\end{proposition}

\proof
The locality condition is equivalent to having \cite{LR}
\be
\hat{\eps}^\pm(\la\otimes\mu^\op,\rho\otimes\nu^\op)
\Psi_{\la,\mu;i} \Psi_{\rho,\nu;j}
= \Psi_{\rho,\nu;j} \Psi_{\la,\mu;i}.
\lableq{lrloc}
Consequently $\Psi_{\la,0;i}$ and $\Psi_{0,\mu;j}$
commute for all $\la,\mu,i,j$, and hence so do
$\hat M _+$ and $\hat M _-^\op$.
\endproof

\subsection{On $\dind$-induction}
\label{Dentire}

Here we show that the entire $R$-$R$ system
can be obtained by some induction procedure similar
to $\a$-induction.

We now define an extension of $\la\otimes\mu^\op\in\End(Q)$
to $R$ in the spirit of $\a$-induction \cite{LR,BE1}
and call it $\dind$-induction:
We put $\dind_{\la,\mu}^\pm|_Q=\la\otimes\mu^\op$ and
\[
\dind_{\la,\mu}^\pm(\Psi_{\rho,\nu;i})=
\hat{\eps}^\pm(\la\otimes\mu^\op,\rho\otimes\nu^\op)^*
\Psi_{\rho,\nu;i}.
\]
We also assume ``locality'', i.e.\
$\hat{\eps}^\pm(\Theta,\Theta)V^2=V^2$.
By the same arguments as for $\a$-induction, this implies
\[
\Hom(\dind^\pm_{\la,\mu},\dind^\pm_{\rho,\nu})
=\Hom(\iota\circ(\la\otimes\mu^\op),
\iota\circ(\rho\otimes\nu^\op))
\]
since \erf{lrloc} implies the opposite inclusion
of the obvious one obtained from the extension property
of $\dind$-induction.
In particular
$\langle\id,\dind_{\la,\mu}^\pm\rangle=Z_{\la,\mu}$.
If we have a subfactor $N\subset M$ allowing for
$\a$-induction then we thus find
$\langle\dind^\pm_{\la,\id},\dind^\pm_{\id,\mu}\rangle
=\langle\a^+_\la,\a^-_{\co\mu}\rangle$.
Therefore we expect $\dind^+_{\la,\mu}$ and
$\dind^-_{\la,\mu}$ to produce systems isomorphic
to $\MXM$. We will show soon that the two sector systems
obtained from $\dind^+$- and $\dind^-$-induction only
intersect on the identity sector.

Let $\RXR$ denote a system of $R$-$R$ morphisms,
i.e.\ a choice of representative morphisms of
irreducible subsectors of some
$[\iota_\mathrm{LR}(\la\otimes\mu^\op)\co\iota _\mathrm{LR}]$
(and containing the identity morphism).
Let $\RXRpm\subset\RXRi\subset\RXR$ denote the
subsystems arising from chiral inductions $\dind^\pm$
and the mixed induction (with notation analogous
to \cite{BEK1,BEK2,BE4,BEK3}).

\begin{lemma}
\label{decouple}
Provided that the braiding on $\NXN$ is non-degenerate,
we have
\be
\langle\dind_{\la,\mu}^+,\dind_{\rho,\nu}^-\rangle
=\langle\dind_{\la,\mu}^+,\id\rangle
\langle\id,\dind_{\rho,\nu}^-\rangle.
\ee
In other words, writing
$\Xi_{\la,\mu;\rho,\nu}=
\langle\dind_{\la,\mu}^+,\dind_{\rho,\nu}^-\rangle$
then gives
$\Xi_{\la,\mu;\rho,\nu}=Z_{\la,\mu}Z_{\rho,\nu}$.
\end{lemma}

\proof
We define a matrix
\[
R_{\la,\mu;\rho,\nu;\Om}^{\Om'} =
\langle \Om \dind^+_{\la,\mu} \dind^-_{\rho,\nu}, \Om' \rangle
\]
for $\Om,{\Om'}\in\RXRi$. Since we assume non-degeneracy of the
braiding, the characters of the $N$-$N$ will be exhausted
by the statistics characters. Thanks to the homomorphism
property of $\dind$-induction with respect to both entries,
we obtain by the usual arguments
\[
R_{\la,\mu;\rho,\nu;\Om}^{\Om'} = \sum_{i=1}^{\# \RXRi}
\frac{S_{\la,\Phi^+_1(i)}}{S_{0,\Phi^+_1(i)}}
\frac{S_{\mu,\Phi^+_2(i)}}{S_{0,\Phi^+_1(i)}}
\frac{S_{\rho,\Phi^-_1(i)}}{S_{0,\Phi^-_1(i)}}
\frac{S_{\nu,\Phi^-_2(i)}}{S_{0,\Phi^-_1(i)}}
\xi^i_\Om (\xi^i_{\Om'})^*
\]
with maps $\Phi^\pm_j:\{1,2,...,\#\RXRi\}\rightarrow \NXN$,
and where the $\xi^i$'s are orthonormal vectors
(cf.\ \cite[p.\ 92]{BE2} or \cite[p.\ 203]{BE3}).
Note that the statistical dimension of $\la\otimes\mu^\op$
is given by $\om_\la/\om_\mu$. By the same arguments as
for $\a$-induction (cf.\ \cite[Lemma 3.10]{BE3}) we find
\[
\Xi_{\la,\mu;\rho,\nu} = \frac{\om_\la \om_\nu}{\om_\mu \om_\rho}
\Xi_{\la,\mu;\rho,\nu} \,.
\]
Note that
$\Xi_{\la,\mu;\rho,\nu}=R_{\la,\mu;\co\rho,\co\nu;0}^0$.
Hence we can compute
\[
\begin{array}{l}
\displaystyle\sum_{\la,\mu,\rho,\nu\in\NXN} S_{0,\la} S_{0,\mu}
\Xi_{\la,\mu;\rho,\nu} S_{\rho,0} S_{\nu,0}
= \displaystyle\sum_{\la,\mu,\rho,\nu\in\NXN} S_{0,\la} S_{0,\mu}
\om_\la \om_\mu^* \Xi_{\la,\mu;\rho,\nu}
\om_\rho^* \om_\nu S_{\rho,0} S_{\nu,0} \\[.4em]
= \displaystyle\sum_{\la,\mu,\rho,\nu\in\NXN}
\sum_{i=1}^{\# \RXRi}
S_{0,\la} S_{0,\mu} \om_\la \om_\mu^*
\frac{S_{\la,\Phi^+_1(i)}}{S_{0,\Phi^+_1(i)}}
\frac{S_{\mu,\Phi^+_2(i)}}{S_{0,\Phi^+_1(i)}}
\frac{S_{\co\rho,\Phi^-_1(i)}}{S_{0,\Phi^-_1(i)}}
\frac{S_{\co\nu,\Phi^-_2(i)}}{S_{0,\Phi^-_1(i)}}
| \xi^i_0 |^2 \om_\rho^* \om_\nu
S_{\rho,0} S_{\nu,0} \\[.4em]
= \displaystyle\sum_{i=1}^{\# \RXRi}
\om_{\Phi^+_1(i)} \om^*_{\Phi^+_2(i)} \om^*_{\Phi^-_1(i)}
\om_{\Phi^-_2(i)} | \xi^i_0 |^2 \le \|\xi_0\|^2=1 ,
\end{array}
\]
where we used the modular relation $STS=T^*ST^*$.
On the other hand we clearly have
$\langle\dind_{\la,\mu}^+,\dind_{\rho,\nu}^-\rangle
\ge \langle\dind_{\la,\mu}^+,\id\rangle
\langle\id,\dind_{\rho,\nu}^-\rangle$, i.e.\
$\Xi_{\la,\mu;\rho,\nu}\ge Z_{\la,\mu} Z_{\rho,\nu}$.
Hence
\[
\displaystyle\sum_{\la,\mu,\rho,\nu\in\NXN} S_{0,\la} S_{0,\mu}
\Xi_{\la,\mu;\rho,\nu} S_{\rho,0} S_{\nu,0} \ge
\displaystyle\sum_{\la,\mu,\rho,\nu\in\NXN} S_{0,\la} S_{0,\mu}
Z_{\la,\mu} Z_{\rho,\nu} S_{\rho,0} S_{\nu,0} =1
\]
by modular invariance. Therefore we conclude that
$\Xi_{\la,\mu;\rho,\nu}=Z_{\la,\mu}Z_{\rho,\nu}$.
\endproof

Next we define a vector $\vec{v}$ by putting
\[
v_\Om = \sum_{\la,\mu,\rho,\nu\in\NXN}
d_\la d_\mu d_\rho d_\nu \langle \Om,
\dind^+_{\la,\mu} \dind^-_{\rho,\nu} \rangle ,
\qquad \Om \in \RXRi.
\]
By an analogous calculation as in the proof of
\cite[Prop.\ 3.1]{BE4} we find
$R_{\la,\mu;\rho,\nu} \vec{v} = d_\la d_\mu d_\rho d_\nu \vec{v}$
so that $\vec{v}$ is a scalar multiple of the Perron-Frobenius
eigenvector $\vec{d}$ (with entries $d_\Om$), i.e.\
$\vec{v}=\zeta\vec{d}$, $\zeta\in\bbR$. By $d_0=1$ we notice
$\zeta=v_0$ which can be calculated as
\[
v_0 = \sum_{\la,\mu,\rho,\nu\in\NXN}
d_\la d_\mu d_\rho d_\nu \langle \id,
\dind^+_{\la,\mu} \dind^-_{\rho,\nu} \rangle
= \sum_{\la,\mu,\rho,\nu\in\NXN}
d_\la Z_{\la,\mu} d_\mu
d_\rho Z_{\rho,\nu} d_\nu = w^2
\]
by Lemma \ref{decouple}. (As usual, $w=\sum_\la d_\la^2$
denotes the global index of $\NXN$.)
Thus $v_\Om=w^2 d_\Om$ for all $\Om\in\RXRi$.
On the other hand we know
\[
d_\la d_\mu d_\rho d_\nu = \sum_{\Om\in\RXRi}
\langle \Om, \dind^+_{\la,\mu} \dind^-_{\rho,\nu} \rangle d_\Om
\]
so that
\[
\begin{array}{ll}
w^4 &= \displaystyle\sum_{\la,\mu,\rho,\nu\in\NXN}
d_\la^2 d_\mu^2 d_\rho^2 d_\nu^2
=\sum_{\la,\mu,\rho,\nu\in\NXN} \sum_{\Om\in\RXRi}
d_\la d_\mu d_\rho d_\nu \langle \Om, \dind^+_{\la,\mu}
\dind^-_{\rho,\nu} \rangle d_\Om \\[.4em]
&= \displaystyle\sum_{\Om\in\RXRi} v_\Om d_\Om =
w^2 \sum_{\Om\in\RXRi} d_\Om^2 = w^2 w_\dind
\end{array}
\]
with $w_\dind = \sum_{\Om\in\RXRi} d_\Om^2$.
Thus we have obtained
\begin{corollary}
We have $w_\dind=w^2$, and consequently each $R$-$R$ sector
is realized by $\dind$-induction, i.e.\ $\RXR=\RXRi$.
\end{corollary}
Let $\Om_\pm\in\RXRpm$. Note that then the product
$\Om_+\Om_-$ is irreducible as
\[
\langle\Om_+\Om_-,\Om_+\Om_-\rangle =
\langle\Om_+{\co\Om}_+,\Om_-{\co\Om}_-\rangle =1
\]
since $\RXRp\cap\RXRm=\{\id\}$ by Lemma \ref{decouple}.
Hence we have obtained
\begin{corollary}
The system of $R$-$R$ morphisms is exhausted by
the products $\Om_+\Om_-$ with $\Om_\pm\in\RXRpm$.
\end{corollary}

\section{Finite groups and the McKay correspondence}
\labl{McKay}

\subsection{Background}

The quantum doubles of finite groups provide modular data
whose modular invariants can be precisely analysed \cite{CGR, EP}.
Finite groups are also of course
significant for modular invariants, particularly $\SUz$ modular invariants,
by the ADE classification of finite subgroups of $\SUz$. If $G$ is
a finite subgroup of $\SUz$, then the McKay graph of $\hat G$
is an affine ADE graph. It is  a challenge to relate the ADE
classsification of modular invariants and corresponding subfactor theory
with the classification of the ADE finite subgroups and the corresponding
theory of the resolution of quotient singularities. Here are some modest
contributions in this direction.
See also announcements of \cite{O, O3}.

\subsection{Kostant's result}
\labl{Kostant}

As is well-known, the McKay correspondence associates
the finite subgroups of $\SUz$ with affine Dynkin diagrams.
Using ordinary Coxeter graphs, the correspondence is
as in Table \ref{ADEtable}.
\begin{table}[htb]
\caption{The McKay correspondence}
\begin{center}
  \begin{tabular}{|c|c|c|} \hline &&\\[-.9em]
Dynkin diagram & Subgroup $G\subset\SUz$ & Order of $G$
 \\[-.9em]&&\\ \hline\hline &&\\[-.9em]
$\rmA_\ell$ & $\bbZ_{\ell+1}$ & $\ell+1$ \\ &&\\[-.9em]
$\rmD_\ell$ & $BD_\ell=Q_{\ell-2}$ & $4\ell-8$ \\ &&\\[-.9em]
$\rmE_6$ & $BT=BA_4$ & $24$ \\ &&\\[-.9em]
$\rmE_7$ & $BO=BS_4$ & $48$ \\ &&\\[-.9em]
$\rmE_8$ & $BI=BA_5$ & $120$ \\[-.8em]&&\\
\hline \multicolumn3c {} \\[.05em] \end{tabular}
\end{center}
\label{ADEtable}
\end{table}
Now pick one of these finite subgroups $G\subset\SUz$.
For $\gamma\in\hat G$ (the group dual) we define
power series
\[
f_\gamma (q) = \sum_{j=0}^\infty n^\gamma_j q^j,
\]
where the numbers $n^\gamma_j$ is the multiplicity
of $\gamma$ in the restricted $(j+1)$-dimensional
representation $D_j$ of $\SUz$,
i.e.\ $n^\gamma_j=\langle\res \SUz G {D_j},\gamma\rangle$.
(Note that $n^\gamma_j \le j+1$ so that $f_\gamma$
converges e.g.\ inside the unit disc of $\bbC$.)
It was proven by Kostant \cite{Kos} that
\[
f_\gamma(q) = \frac{p_\gamma(q)}{(1-q^r)(1-q^s)}
\]
with $r$ and $s$ positive integers such that $r+s=h+2$
and $rs=\# G$ where $h$ is the dual Coxeter number of
the Dynkin diagram associated to the group $G$
(with order $\# G$), and the important thing is now
that $p_\gamma(q)$ is a polynomial.
Labelling the trivial representation in $\hat G$ by ``$\ast$''
then we have $p_\ast(q)=1+q^h$ and for the other representations
we have
\[
p_\gamma(q) = \sum_{i=0}^{h-2} c^\gamma_i q^{i+1},
\]
with positive integers $c^\gamma_i$.
But let us now recall that the affine Dynkin diagram
$\hat\Gamma$ associated to $G$ has adjacency matrix
\[
\hat\Gamma _{\gamma,\gamma'} = \langle \gamma',
\gamma\otimes\gamma_1 \rangle , \qquad \gamma,\gamma'\in\hat G,
\]
and where $\gamma_1=\res\SUz G{D_1}$.
Recall the tensor product rules
$D_j\otimes D_1=D_{j+1}\oplus D_{j-1}$ for $j\ge1$
and certainly $D_0\otimes D_1=D_1$.
Thus we can compute
\[
\begin{array}{ll}
\sum_{\gamma'\in\hat G} \hat\Gamma _{\gamma,\gamma'} f_{\gamma'}(q)
&= \sum_{\gamma'\in\hat G} \sum_{j=0}^\infty
\langle\res\SUz G{D_j}, \gamma'\rangle \langle \gamma',
\gamma\otimes\gamma_1\rangle q^j \\[.4em]
&= \sum_{j=0}^\infty \langle\res\SUz G{D_j},
\gamma\otimes\gamma_1\rangle q^j \\[.4em]
&= \sum_{j=0}^\infty \langle\res\SUz G{D_1\otimes D_j},
\gamma\rangle q^j \\[.4em]
&= \sum_{j=0}^\infty n^\gamma_{j+1}q^j + \sum_{j=1}^\infty
n^\gamma_{j-1} q^j \\[.4em]
&=q^{-1}(f_\gamma(q)-n^\gamma_0)+q f_\gamma(q),
\end{array}
\]
so that we obtain
\be
\sum_{\gamma'\in\hat G} \hat\Gamma _{\gamma,\gamma'} f_{\gamma'}(q)
=(q+q^{-1})f_\gamma(q) - q^{-1}\delta_{\gamma,\ast},
\lableq{must}
due to $n^\gamma_0=\delta_{\gamma,\ast}$,
and in fact \erf{must} determines the $f_\gamma$'s
completely \cite{Itz}.
(For $\SUd$, presumably  the polynomials would have two variables,
i.e.\ something like
\[
f_\gamma (q_1,q_2) = \sum_{\la_1=0}^\infty \sum_{\la_2=0}^\infty
n^\gamma_{(\la_1,\la_2)} q_1^{\la_1} q_2^{\la_2}
\]
with $n^\gamma_{(\la_1,\la_2)}
=\langle\res \SUd G {D_{(\la_1,\la_2)}},\gamma\rangle$
for $G\subset\SUd$, $\gamma\in\hat G$, and
$D_{(\la_1,\la_2)}$ is the representation of $\SUd$ with
``Dynkin labels'' $\la_1,\la_2\ge0$

\subsection{Relating Kostant's result to subfactors}

Now let us turn to our subfactors realizing the ADE modular
invariants. Recall that for an ordinary ADE Dynkin diagram
$\Gamma$ its affine extension $\hat\Gamma$ is obtained by
adding one vertex ``$\ast$'' and joining it by an edge to a
certain vertex of $\Gamma$. By the argument of
changing the $\iota$-vertex  \cite{E1}  we may and do assume that
our subfactor $N\subset M$ realizing the $\Gamma$ modular
invariant has the $\iota$-vertex exactly on the vertex
which would join the extension vertex ``$\ast$''.
It turns out that we then have $\langle\canr,\la_k\rangle=1$
for all the DE cases (for A we will need =2 and a reducible
subfactor, what we do  discuss here in detail).
This gives a natural bijection between (equivalence classes of)
non-trivial irreducible representations of $G$ and $M$-$N$
sectors. We will denote the $M$-$N$ morphism associated to
$\gamma\in\hat G\setminus\{\ast\}$ by $\co a _\gamma$.
Note that $\iota=\co a _{\gamma_1}$.
Now let us define the following polynomials labelled by
$\gamma\in\hat G$:
\begin{equation}
\begin{array}{rl}
p_\ast(q) &= 1 + q^{k+2} ,\\[.4em]
p_\gamma (q) &= \sum_{j=0}^k \langle \co a _\gamma,
 \iota\la_j \rangle q^{j+1}, \\[.4em]
\Omega(q) &= (1+q^2) p_\ast(q) - qp_{\gamma_1}(q).
\end{array}
\lableq{polys}
We now claim
\[
\sum_{\gamma'\in\hat G} \hat\Gamma _{\gamma,\gamma'} p_{\gamma'}(q)
=(q+q^{-1})p_\gamma(q) - q^{-1}\delta_{\gamma,\ast}\Omega(q),
\]
Let's compute the left-hand side for the three cases:
(1) $\gamma=\ast$,
(2) $\gamma=\gamma_1$,
and (3) otherwise.
For (1) we clearly obtain $p_{\gamma_1}$.
For (2) we obtain
$p_\ast+\sum_{\gamma'}\langle\co a _\gamma\la_1,
\co a _{\gamma'}\rangle p_{\gamma'}$,
and for (3) we obtain
$\sum_{\gamma'}\langle\co a _\gamma\la_1,
\co a _{\gamma'}\rangle p_{\gamma'}$.
Now note that
\[
\begin{array}{rl}
\sum_{\gamma'\in\hat G}\langle\co a _\gamma\la_1,
\co a _{\gamma'}\rangle p_{\gamma'}(q)
&= \sum_{j=0}^k \langle\co a _\gamma\la_1,
\iota\la_j \rangle q^{j+1} \\[.4em]
&= \sum_{j=0}^{k-1} \langle\co a _\gamma,
\iota\la_{j+1} \rangle q^{j+1} +
\sum_{j=1}^k \langle\co a _\gamma,
\iota\la_{j-1} \rangle q^{j+1} \\[.4em]
&=(q+q^{-1})p_\gamma(q) - \langle\co a _\gamma,
\iota \rangle  - \langle\co a _\gamma,
\iota\la_k \rangle q^{k+2} \\[.4em]
&= (q+q^{-1})p_\gamma(q) -(1+q^{k+2})\delta_{\gamma,\gamma_1},
\end{array}
\]
since $[\iota]=[\iota\la_k]$. Therefore the left-hand side
gives $p_{\gamma_1}(q)$ for (1), i.e.\ $\gamma=\ast$, and
$(q+q^{-1})p_\gamma(q)$ for both (2) and (3), i.e.\ otherwise.
This equals the right-hand side, so the claim is correct.
We conclude that putting $P_\gamma(q)=p_\gamma(q)/\Omega(q)$
in fact solves \erf{must}, and thus we can identify
$f_\gamma=P_\gamma$ and in turn
$c^\gamma_i=\langle\co a _\gamma,\iota\la_j\rangle$.
Thus we have a fairly direct interpretation of the
$c^\gamma_i$'s through putting the $\iota$-vertex next
to the $\ast$-vertex.
But note that even if we locate a different vertex,
say $\bar b$, next to $\ast$, then we will similarly obtain
$c^\gamma_i=\langle\co a _\gamma,\co b\la_j\rangle$.
(Therefore the $c^\gamma_i$'s are found to be certain
entries of our matrices $G_\la$. This coincidence has
also been noticed by Di Francesco and Zuber \cite{DiF}
in their more empirical approach to these matrices.)
Nevertheless the exponents in Kostant's polynomials
appear most directly in the above sketched setting
where $p_{\gamma_1}$ ``is essentially $[\canr]$''.
(Note that $[\canr]=\bigoplus_j c^{\gamma_1}_j[\la_j]$.)

\subsection{A modular invariant from a group}

We now construct modular invariant for a simple case
without $\a$-induction.
Let $\Delta\subset N$ be a braided system, $\Gamma\subset\Delta$
a subsystem and let $\Theta\subset\Gamma$ be the degenerate
subsystem of $\Gamma$. Thanks to Rehren \cite{R1}, and the theory of
Doplicher-Roberts \cite{Dop}, we know that
$\Theta$ is a group dual, i.e.\ $\om_\vartheta=\pm 1$ and
$d_\vartheta\in\bbZ$ for all $\vartheta\in\Theta$.
We assume that it is purely bosonic, i.e.\ that
$\om_\vartheta=+1$ for all $\vartheta\in\Theta$.
Let $w_\Gamma=\sum_{\gamma\in\Gamma}d_\gamma^2$
and put
\[
\langle y^\la,y^\mu\rangle_\Gamma
=\sum_{\gamma\in\Gamma} Y_{\la,\gamma}^* Y_{\mu,\gamma}.
\]
Note that for $\la,\mu\in\Gamma$ we have
\[
\langle y^\la,y^\mu\rangle_\Gamma=\sum_{\gamma\in\Gamma}
Y_{\co\la,\gamma}Y_{\mu,\gamma}=\sum_{\gamma\in\Gamma}
d_\gamma \sum_{\gamma'\in\Gamma} N_{\co\la,\mu}^{\gamma'}
Y_{\gamma,\gamma'} = \sum_{\gamma'\in\Gamma}
N_{\co\la,\mu}^{\gamma'} \langle y^0,y^{\gamma'}\rangle_\Gamma
= w_\Gamma \sum_{\vartheta\in\Theta}
N_{\la,\vartheta}^\mu d_\vartheta
\]
Let's now define a matrix $Z$ with entries
\[
Z_{\la,\mu} = \left\{ \begin{array}{cl}
w_\Gamma^{-1} \langle y^\la,y^\mu\rangle_\Gamma &
\qquad \la,\mu\in\Gamma,\\
0 & \qquad \mbox{otherwise,}
\end{array} \right.
\]
which hence has non-negative integral entries,
and also $Z_{0,0}=1$.
Since $\Theta$ is the subsystem of degenerate sectors
in $\Gamma$ we will have $\om_\la \om_\vartheta /\om_\mu=1$
whenever $N_{\la,\vartheta}^\mu\neq0$, and since
it is purely bosonic we will then therefore have
$\om_\la=\om_\mu$ for $\la,\mu\in\Gamma$.
It follows that $ZT=TZ$. To establish Y-invariance,
we need another assumption. We assume that
$\langle y^0,y^\mu\rangle_\Gamma=0$
whenever $\la\notin\Gamma$.
Note that this implies for any $\gamma\in\Gamma$
that $\langle y^\la,y^\gamma\rangle_\Gamma=0$ whenever
$\la\notin\Gamma$ since then
\[
\begin{array}{rl}
\langle y^\la,y^\gamma\rangle_\Gamma &= \sum_{\gamma'\in\Gamma}
Y_{\co\la,\gamma'}Y_{\gamma,\gamma'} = \sum_{\gamma'\in\Gamma}
d_{\gamma'} \sum_{\nu\in\Delta} N_{\co\la,\gamma}^\nu
Y_{\gamma',\nu} \\[.4em]
&= \sum_{\nu\in\Delta} N_{\co\la,\gamma}^\nu
\langle y^0,y^\nu\rangle_\Gamma =\sum_{\gamma'\in\Gamma}
N_{\co\la,\gamma}^{\gamma'} \langle y^0,
y^{\gamma'} \rangle_\Gamma=0
\end{array}
\]
since then $N_{\co\la,\gamma}^{\gamma'}=0$.
Now Y-invariance can be easily seen: Writing $Y$ and $Z$
as $2\times2$ blocks according to labels in $\Gamma$ and
not in $\Delta\setminus\Gamma$ we need to check that
\[
\sum_{\gamma\in\Gamma} Y_{\la,\gamma} Z_{\gamma,\mu}=
\sum_{\gamma\in\Gamma} Z_{\la,\gamma} Y_{\gamma,\mu}
\]
for $\la,\mu\in\Gamma$, and that
\[
\sum_{\gamma\in\Gamma} Y_{\la,\gamma} Z_{\gamma,\mu}=0
\]
whenever $\la\notin\Gamma,\mu\in\Gamma$.
Using our assumption, both is now obvious if we write
\[
\sum_{\gamma\in\Gamma} Y_{\la,\gamma} Z_{\gamma,\mu}
= w_\Gamma^{-1} \sum_{\gamma,\rho\in\Gamma}
Y_{\la,\gamma} Y_{\gamma,\rho}^* Y_{\rho,\mu}
= w_\Gamma^{-1} \sum_{\rho\in\Gamma}
\langle y^\la, y^\rho \rangle_\Gamma Y_{\rho,\mu}
\]
for $\mu\in\Gamma$.

\section{How to prove $YZ=ZY$ without the double triangle algebra}
\labl{algebraic}

Here we show  how one can prove Y-invariance
of the coupling matrix ($\Om$-invariance is established
in \cite[Lemma 6.1]{BE4}) without reference to the
double triangle algebra (DTA). However, we can only complete
the argument covering the non-degenerate case and
the case of $Z$ being a permutation.

We have the nimrep $\{\Gamma_\lambda : \lambda \in \NXN \}$ of multiplication by
$\alpha^{\pm}_\lambda$ on $\MXMpm$. Taking quantum dimension, since $d_\lambda =
d_{{\alpha}_\lambda}$, we must have that $N_\lambda[d_\beta]_\beta
= d_\lambda[d_\beta]_\beta$, using our basis $\beta$ of $\MXMpm$. On the other hand,
$N_\lambda\alpha^{\pm}_\mu = \alpha^{\pm}_{\lambda\mu}$, and again we
should get
a vector [$d_\mu$] with eigenvalue $d_\lambda$. To relate this to the
basis $[\beta]$, we just proceed as in the \ Subsection \ref{Dentire}

First we recall from the proof of \cite[Prop.\ 3.1]{BEK2}
that by defining vectors $\vec{v}^\pm$ with components
labelled by $\beta\in\MXMpm$,
\[
v^\pm_\beta = \sum_{\la\in\NXN} d_\la
\langle\beta,\a^\pm_\la\rangle,
\qquad\beta\in\MXMpm,
\]
a simple Perron-Frobenius argument shows that the vectors
$\vec{v}^\pm$ are proportional to the dimension vectors
$\vec{d}^\pm$ (with entries $d_\beta$), and since
$d_\la=\sum_\beta \langle\beta,\a^\pm_\la\rangle d_\beta$
we find by summing $\sum_\beta d_\beta v^\pm_\beta$ that
$\vec{v}^\pm=\vec{d}^\pm w/w_\pm$. Comparing the
$0$-entries yields then
\[
\frac w{w_+} = \sum_{\la\in\NXN} d_\la Z_{\la,0}, \qquad
\frac w{w_-} = \sum_{\la\in\NXN} Z_{0,\la} d_\la.
\]
It is obvious that commutativity $[Y,Z]=0$
(which was proven in \cite[Thm.\ 5.7]{BEK1} using DTA methods)
now implies $w_+=w_-$,
but this is what we
want  to avoid here since the aim
is now to prove $[Y,Z]=0$ independently.

By a similar Perron-Frobenius argument,
it was shown in the proof of \cite[Prop.\ 3.1]{BE4}
that the vector $\vec{v}$ with entries
\[
v_\beta = \sum_{\la,\mu\in\NXN}
d_\la d_\mu \langle\beta,\a^+_\la\a^-_\mu\rangle,
\qquad\beta\in\MXMa,
\]
is proportional to the vector $\vec{d}$, now with entries
$d_\beta$ labelled by $\beta\in\MXMa$, and one finds
that $\vec{v} = \vec{d} w^2/w_\a$.
Checking the $0$-entries one thus obtains
\[
\frac{w^2}{w_\a} = \sum_{\la,\mu\in\NXN} d_\la Z_{\la,\mu} d_\mu.
\]
(Extending the vectors $\vec{v}^\pm$ to $\MXMa$ by entries $=0$
and computing $\langle\vec{v}^+,\vec{v}^-\rangle$ also yields
$w_\a w_0=w_+w_-$.)
Now as in the proof of \cite[Lemma 6.2]{BE4} we define
two further vectors $\vec{u}^\pm$ with entries
\[ u^\pm_\beta=\sum_{\la\in\NXN} \om_\la d_\la
\lan\beta,\a^\pm_\la\ran \,, \qquad \beta\in\MXMa \,.\]
(Or labelled by $\MXM$, it does not  matter.)
We now want to establish a very important identity,
namely
\be
u^\pm_\beta
= \left\{ \begin{array}{c@{\qquad:\qquad}l}
 \om_\tau d_\tau w/w_\pm & \beta=\tau\in\MXMo \\[.6em]
 0 & \mbox{otherwise.} \end{array} \right.
\label{inphsum}
\ee
To prove this, we first compute
\[ \begin{array}{ll}
\| \vec{u}^+ \|^2 &= \sum_\beta \sum_{\la,\nu} \om_\la \om_\nu^{-1}
 d_\la d_\nu \lan \a^+_\nu,\beta \ran \lan \beta, \a^+_\la \ran
= \sum_{\la,\nu} \om_\la \om_\nu^{-1}
 d_\la d_\nu \lan \a^+_{\co\la} \a^+_\nu, \id \ran \\[.4em]
&= \sum_{\la,\mu,\nu} \om_\la \om_\nu^{-1}
 d_\la d_\nu N_{\co\la,\nu}^\mu \lan \a^+_\mu, \id \ran
= \sum_{\la,\mu,\nu} \om_\la \om_\nu^{-1}
 d_\la d_\nu N_{\la,\mu}^\nu \om_\mu Z_{\mu,0} \\[.4em]
&= \sum_{\la,\mu} Y_{0,\la} Y_{\la,\mu} Z_{\mu,0} \,,
\end{array} \]
where we used that $Z_{\mu,0}=\om_\mu Z_{\mu,0}$
by \cite[Lemma 6.1]{BE4}. Similarly we obtain
$\|\vec{u}^-\|^2=\sum_{\la,\mu}Y_{0,\la}Y_{\la,\mu}Z_{0,\mu}$.
Now by Rehren's argument \cite{R1} we have
\[ \sum_{\mu\in\NXN} Y_{\la,\mu} Y_{\mu,0} = \left\{
\begin{array}{cl} w d_\la & \qquad \la\in\NXNd \\
0 & \qquad \la\notin\NXNd \end{array} \right. . \]
Since clearly $Z_{\la,0}=Z_{0,\la}$ as
$\a^+_\la=\a^-_\la$ for $\la\in\NXNd$,
we obtain
\be
\|\vec{u}^+\|^2 = \|\vec{u}^-\|^2
= w \sum_{\la\in\NXNd} d_\la Z_{\la,0} .
\label{++--}
\ee
On the other hand, we can compute the inner product
\be
\begin{array}{rl}
\lan \vec{u}^+,\vec{u}^- \ran
&= \sum_{\la,\mu} \sum_\beta \om_\la^{-1} d_\la d_\mu \om_\mu
\lan \a^+_\la,\beta\ran \lan\beta,\a^-_\mu \ran
= \sum_{\la,\mu} \om_\la^{-1} d_\la d_\mu \om_\mu
 Z_{\la,\mu} \\[.4em]
&= \sum_{\la,\mu} d_\la Z_{\la,\mu} d_\mu
= \displaystyle\frac{w^2}{w_\a}.
\end{array}
\label{+-}
\ee
Now let us consider the case that the braiding is
non-degenerate or that $Z$ is a permutation matrix
(i.e.\ $Z_{\la,0}=\delta_{\la,0}$).
Then \erf{++--} yields $\|\vec{u}^+\|^2=\|\vec{u}^-\|^2=w$
whereas \erf{+-} yields for the inner product
$\lan \vec{u}^+,\vec{u}^- \ran=w^2/w_\a \ge w$
(since always $w_\a\le w$), so that we can
conclude $\vec{u}^+=\vec{u}^-$.
Since obviously $u^\pm_\beta=0$ whenever $\beta\notin\MXMpm$
this implies $u^\pm_\beta=0$ whenever $\beta\notin\MXMo$.
Now for $\beta=\tau\in\MXMo$ we find
by \cite[Lemma 6.1]{BE4} and the above derived
$\vec{v}^\pm=\vec{d}^\pm w/w_\pm$
\[
u^\pm_\tau = \sum_{\la\in\NXN} \om_\la d_\la
\lan\tau,\a^\pm_\la\ran = \om_\tau \sum_{\la\in\NXN}
d_\la \lan\tau,\a^\pm_\la\ran =
\frac w{w_\pm} d_\tau \om_\tau ,
\]
proving \erf{inphsum}.
 So what is so fascinating
about \erf{inphsum}? The point is, that it fully
equivalent to $[Y,Z]=0$. In fact, using $[Y,Z]$
we proved $\vec{u}^+=\vec{u}^-$ and thus \erf{inphsum}
in \cite[Lemma 6.2]{BE4}.
Conversely, \erf{inphsum} can be used, as in
the proof of \cite[Lemma 6.3]{BE4}, to compute

\[ \begin{array}{rl}
\displaystyle\sum_{\la\in\NXN} b^\pm_{\tau,\la} Y_{\la,\mu}
&= \displaystyle\sum_{\la,\rho\in\NXN} \frac{\om_\la \om_\mu}{\om_\rho}
 N_{\la,\mu}^\rho d_\rho b^\pm_{\tau,\la} \\[2em]
&= \displaystyle\sum_{\la,\rho\in\NXN} \frac{\om_\tau \om_\mu}{\om_\rho}
 N_{\rho,\co\mu}^\la d_\rho \lan\a^\pm_\la, \tau \ran \\[2em]
&= \displaystyle\sum_{\rho\in\NXN}
 \frac{\om_\tau \om_\mu}{\om_\rho} d_\rho
 \lan\a^\pm_\rho, \tau \a^\pm_\mu\ran \\[2em]
&= \displaystyle\sum_{\rho\in\NXN} \sum_{\beta\in\MXMpm}
 \frac{\om_\tau \om_\mu}{\om_\rho} d_\rho \lan\a^\pm_\rho,
 \beta \ran \lan \beta, \tau \a^\pm_\mu\ran \\[2em]
&= \displaystyle\frac w{w_\pm} \sum_{\tau''\in\MXMo}
 \frac{\om_\tau \om_\mu}{\om_{\tau''}}
  d_{\tau''} \lan \tau'',\tau \a^\pm_\mu\ran \\[2em]
&= \displaystyle\frac w{w_\pm} \sum_{\tau',\tau''\in\MXMo}
 \frac{\om_\tau \om_\mu}{\om_{\tau''}} N_{\co\tau,\tau''}^{\tau'}
  d_{\tau''} \lan \tau', \a^\pm_\mu\ran \\[2em]
&= \displaystyle\frac w{w_\pm} \sum_{\tau',\tau''\in\MXMo}
 \frac{\om_\tau \om_{\tau'}}{\om_{\tau''}}
 N_{\tau,\tau'}^{\tau''} d_{\tau''} b^\pm_{\tau',\mu} \\[2em]
&= \displaystyle\frac w{w_\pm} \sum_{\tau'\in\MXMo}
 Y^\ext_{\tau,\tau'} b^\pm_{\tau',\mu} \,,
\end{array} \]
where we actually used (the complex conjugate of) \erf{inphsum}
in the fifth equality. Here we were however a little more careful
than in the proof of \cite[Lemma 6.3]{BE4} as we here
distinguished $w_+$ and $w_-$ because we would not know
about their equality without using \cite[Thm.\ 5.7]{BEK1}.
Writing the chiral branching coefficients as rectangular matrices,
we have now obtained $w Y^\ext b^\pm=w_\pm b^\pm Y$.
With $Z=\tmat {b^+}b^-$ this yields
\[
w_- ZY= w_- \tmat {b^+}b^- Y = w \, \tmat {b^+}Y^\ext b^-
= w_+ Y \tmat {b^+}b^- =w_+ YZ .
\]
Since $Z_{\la,0}=Z_{0,\la}$ for $\la\in\NXNd$,
we obtain in particular
\[ \frac{w_+}{w_-}(YZY)_{0,0}=(ZY^2)_{0,0}
= \sum_{\la\in\NXNd} Z_{0,\la}d_\la
= (Y^2Z)_{0,0} = \frac{w_-}{w_+}(YZY)_{0,0} \,. \]
Now $(YZY)_{0,0}\ge Z_{0,0}=1$, hence $w_+=w_-$
so that we finally obtain $YZ=ZY$ indeed.

To summarize, we have shown
that the following are equivalent:
\begin{itemize}
\item $[Y,Z]=0$,
\item $\sum_{\la\in\NXNd} d_\la Z_{\la,0}\le w/w_\a$, and
\item $\sum_{\la\in\NXN}\om_\la d_\la\lan\beta,\a^\pm_\la\ran=0$
 unless $\beta$ is ambichiral.
\end{itemize}
If $Z_{\la,0}$ vanishes for all degenerate $\la$ except
$\la=\id$ (in particular if the braiding is non-degenerate),
then the second condition is trivially met since
clearly $1\le w/w_\a$.
Unfortunately, for the general case we still need the DTA
to prove the first condition (and hence all).

\vspace{0.2cm}\addtolength{\baselineskip}{-2pt}
\begin{footnotesize}
\noindent{\it Acknowledgement.}
This work was completed during visits to the MSRI programme
on Operator Algebras in 2000-2001.
 I am grateful for the organizers
and MSRI for their invitation and financial support. I would like to thank
Jens B\"ockenhauer for discussions on this work
as well as for our collaboration in 1996--2000 and Karl-Henning Rehren
for his comments on a preliminary version of this manuscript.
\end{footnotesize}
\vspace{0.5cm}



\newcommand\bitem[2]{\bibitem{#1}{#2}}

\def\aam              {Acta Appl.\ Math. }
\def\aip              {Ann.\ Inst.\ H.\ Poincar\'e (Phys.\ Th\'eor.) }
\def\cmp              {Com\-mun.\ Math.\ Phys. }
\def\duke             {Duke Math.\ J. }
\def\ijm              {Intern.\ J. Math. }
\def\jfa              {J.\ Funct.\ Anal. }
\def\jmp              {J.\ Math.\ Phys. }
\def\lmp              {Lett.\ Math.\ Phys. }
\def\rmp              {Rev.\ Math.\ Phys. }
\def\inv              {Invent.\ Math. }
\def\mpl              {Mod.\ Phys.\ Lett. }
\def\nup              {Nucl.\ Phys. }
\def\nupp             {Nucl.\ Phys.\ (Proc.\ Suppl.) }
\def\npbp             {Nucl.\ Phys.\ {\bf B} (Proc.\ Suppl.)}
\def\adma             {Adv.\ Math. }
\def\physa            {Physica {\bf A} }
\def\ijmp             {Int.\ J.\ Mod.\ Phys. }
\def\jp               {J.\ Phys. }
\def\fdp              {Fortschr.\ Phys. }
\def\plb              {Phys.\ Lett.\ {\bf B}}
\def\rims             {Publ.\ RIMS, Kyoto Univ. }


\begin{footnotesize}

\end{footnotesize}
\end{document}